%detconj.tex: A  Conjectured Explicit Determinant  Evaluation Whose Proof 
%%a Plain TeX file by Doron Zeilberger (x pages)

%begin macros

\baselineskip=14pt
\parskip=10pt

\font\eightrm=cmr8 
\font\eighttt=cmtt8
\magnification=\magstephalf

\def\1{{\overline{1}}}
\def\2{{\overline{2}}}
\parindent=0pt
\overfullrule=0in

\def\frac#1#2{{#1 \over #2}}
%\headline={\rm  \ifodd\pageno  \RightHead  \else  \LeftHead  \fi}
%\def\RightHead{\centerline{
%Title
%}}
%\def\LeftHead{ \centerline{Doron Zeilberger}}
%end macros
\centerline
{
\bf
A  Conjectured Explicit Determinant  Evaluation Whose Proof 
}
\centerline
{
\bf Would Make Me Happy  (and the OEIS richer)
}

\rm
\bigskip
\centerline{ {\it
Doron 
ZEILBERGER}\footnote{$^1$}
{\eightrm  \raggedright
Department of Mathematics, Rutgers University (New Brunswick),
Hill Center-Busch Campus, 110 Frelinghuysen Rd., Piscataway,
NJ 08854-8019, USA.
%\break
{\eighttt zeilberg  at math dot rutgers dot edu} ,
\hfill \break
{\eighttt http://www.math.rutgers.edu/\~{}zeilberg/} .
Jan. 7, 2014.
downloadable from Zeilberger's website.
Supported in part by the NSF.
}
}

{\bf Abstract}: I conjecture a certain explicit determinant evaluation, whose proof would imply the solution
of  certain enumeration problem that I have been working on, and that I find interesting.
I am pledging \$500 to the OEIS Foundation (in honor of the prover!) for a proof,
and \$50 (in honor of the disprover or his or her computer) for a disproof, as well as
(in the affirmative case only) a co-authorship in a good enumeration paper, that would  immediately
bequest a Zeilberger-number $1$, an Erd\"os number $\leq 3$, an Einstein number $\leq 4$ and
numerous other prestigious numbers.

In order to complete the proof of a certain enumeration problem that I have been working on
for the last few weeks, I need a proof of the following conjecture.

Let $d$ be a positive integer, and
Let $M=M(d)$ be the following $2d \times 2d$ matrix with entries in  $\{-1,0,1\}$.
For $1\leq a \leq 2d$ and $1 \leq b \leq d$, 
$$
M_{a,2b-1}=
 \cases{
1 & if $a=2b$ ;\cr
 -1      &  if  $a=3b+1$;\cr
0 . & otherwise
}   
$$
$$
M_{a,2b}=
 \cases{
1 & if $a=2b-1$ ;\cr
 -1      &  if  $a=b-1$;\cr
0 . & otherwise
}   
$$

{\bf Conjecture}:  For every positive integer $d$, the following is true:
$$
\det M(d)=(-1)^d \quad .
$$

{\bf Comments}: 

{\bf 1.} This conjecture came up in my current work in enumerative combinatorics.
Shalosh B. EKhad kindly verified it for $d \leq 200$.
I have no idea how hard it is, and it is possibly not that hard, but right now
I am busy with other problems. I believe that the powerful and versatile techniques
of Karattenthaler[K1][K2] may be applicable, and possibly the
computer-assisted approach described in [Z] and already 
nicely exploited in [KKZ] and [KT].

{\bf 2.} The Short Maple code in:
{\tt http://www.math.rutgers.edu/\~{}zeilberg/tokhniot/DetConj}
defines the matrix $M(d)$ (procedure {\tt M(d)}) and procedure
{\tt C(N)} verifies it empirically for all $d \leq n$.
So far {\tt C(200);} returned {\tt true}.

{\bf 3.} I am offering to donate \$500 to the OEIS Foundation for a proof and \$50 for a disproof,
with an explicit statement that the donation is in honor of the prover (or disprover).

{\bf References}

[K1] Christian Krattenthaler, {\it Advanced Determinant Calculus},
S\'em. Lothar. Comb. {\bf 42} (1999), B42q.
(``The Andrews Festschrift'', D. Foata and G.-N. Han
(eds.))
\hfill\break
{\tt http://www.mat.univie.ac.at/~kratt/artikel/detsurv.html}

[K2] Christian Krattenthaler, {\it Advanced Determinant Calculus: a complement},
Linear Algebra Appl. {\bf 411} (2005), 68-166
\hfill\break
{\tt http://www.mat.univie.ac.at/~kratt/artikel/detcomp.html}

[KKZ] Christoph Koutschan, Manuel Kauers, and Doron Zeilberger , {\it A Proof Of George Andrews' and David Robbins' q-TSPP Conjecture},
Proceedings of the National Academy of Science, {\bf 108\#6} (Feb. 8, 2011), 2196-2199.

[KT] Christoph Koutschan and Thotsaporn Thanatipanonda,
{\it Advanced Computer Algebra for Determinants.},
Annals of Combinatorics, to appear.
\hfill\break
{\tt http://www.risc.jku.at/people/ckoutsch/det/}

[Z] Doron Zeilberger, {\it The HOLONOMIC ANSATZ II: Automatic DISCOVERY(!) and PROOF(!!) of Holonomic Determinant Evaluations},
Annals of Combinatorics {\bf 11} (2007), 241-247
\hfill\break
{\tt http://www.math.rutgers.edu/\~{}zeilberg/mamarim/mamarimhtml/ansatzII.html}

\end